\documentclass[a4paper]{article}
\usepackage{latexsym, amsmath, amssymb}
\usepackage{graphicx,epsfig}

\usepackage{amsfonts}
\usepackage{amscd}
\usepackage{amsbsy}
\usepackage{bm}
\usepackage{color}

\newcommand{\bE}{\mathbf{E}}

\newcommand{\bH}{\mathbf{H}}

\newcommand{\Bx}{\mathbf{x}}
\newcommand{\By}{\mathbf{y}}

\newcommand{\p}{\partial}

\newcommand{\Gl}{\lambda}

\newcommand{\Gs}{\sigma}

\newtheorem{thm}{Theorem}[section]

\newtheorem{lem}[thm]{Lemma}

\newtheorem{defn}[thm]{Definition}

\numberwithin{equation}{section}
\newcounter{saveeqn}

\def\nm{\noalign{\medskip}}
\newcommand{\eqnref}[1]{(\ref {#1})}

\newcommand{\Bb}{\mathbf{b}}

\newcommand{\Kcal}{\mathcal{K}}

\newcommand{\Scal}{\mathcal{S}}
\newcommand{\Mcal}{\mathcal{M}}

\newcommand{\ds}{\displaystyle}

\newcommand{\pf}{\medskip \noindent {\sl Proof}. \ }
\newcommand{\qed}{\hfill $\Box$ \medskip}
\newcommand{\RR}{\mathbb{R}}

\newcommand{\beq}{\begin{equation}}
\newcommand{\eeq}{\end{equation}}

\DeclareMathAlphabet{\itbf}{OML}{cmm}{b}{it}

\def\by{{{\mathbf{y}}}}
\def\bx{{{\mathbf{x}}}}

\def\bz{{{\itbf z}}}

\begin{document}
\title{{\itshape Plasmon resonance and heat generation model in nanostructures}
\thanks{Fang is supported by NSF grants  No. NSFC70921001 and No. 71210003. Deng and Li are supported by NSF grants No. NSFC11301040,}}
\author{Xiaoping Fang \thanks{Postdoctoral, Management Science and Engineering Postdoctoral Mobile Station, School of Business; School of Mathematics and Statistics, Central South University, Changsha, Hunan 410083, P. R. China. Email: fxpmath@csu.edu.cn} \quad
    Youjun Deng \thanks{Corresponding author. School of Mathematics and Statistics,
   Central South University, Changsha, Hunan 410083, P. R. China. Email: youjundeng@csu.edu.cn, dengyijun\_001@163.com}  \quad
    Jing Li \thanks{Department of Mathematics, Changsha University of Science and Technology, Changsha, Hunan 410004, P.R.China. Email: lijingnew@126.com}
}
\date{}
\maketitle

\begin{abstract}
In this paper, we investigate the photothermal effects of the plasmon resonance. Metal nanoparticles
efficiently generate heat in the presence of electromagnetic radiation. The process is strongly enhanced
when a fixed frequency of the incident wave illuminate on nanoprticles such that plasmon resonance happen.
We shall introduce the
electromagnetic radiation model and show exactly how and when the plasmon resonance happen. We then
construct the heat generation and transfer model and derive the heat effect induced by plasmon resonance.
Finally, we consider the heat generation under plasmon resonance in a concentric nanoshell structure.\\
\noindent {\footnotesize {\bf Mathematics subject classification
(MSC2000):} 35Q60, 35C20}
\end{abstract}

\section{Introduction}
There has been a great deal of interest in recent years in the study of physical property - heat
generation by Nanoparticles (NPs) under optical illumination. The optical properties of NPs, including both
semiconductor \cite{MNB1993} and metal nanocrystals \cite{ESMLM1997}, have been studied intensively.
When the NPs are illuminated by an incident wave with special frequency, the plasmon resonance
can be excited with a strong field enhancement inside and around the NPs. Such effect has a wide range of
potential applications in such as near-field microscopy \cite{BL2002,BLBCL2006}, signal amplification, molecular recognition \cite{HF04},
nano-lithography \cite{Pe99}, ect.  Heat is generated when imaginary part of the property of the NP does not vanish. The heating effect becomes especially strong under the
plasmon resonance conditions when the energy of incident photons is close to the plasmon frequency of the NP. The heat can be used
for melting the surrounding matrix like ice and polymer \cite{GR07,Ka13,WMTH2009}, as well as for cancer diagnosis and therapy \cite{GLHJDW2007,HEQE2006,JS07,ZK2009}.

To mathematically explain the plasmon resonance of NPs and the heat generated by plasmon resonance, there are
two models to be considered. Firstly, the Maxwell system is introduced to explain the process of plasmon resonance.
We shall give the exact mathematical theory on how the plasmon resonance happen
when the NPs is illuminated by electromagnetic wave. We show the asymptotic expansion of electromagnetic fields perturbed by small inclusions
derived in \cite{ADKLMZ, AVV2001}. The first order expansion is quite important to see the innate of the inclusion.
It can be used for the size estimation of the small inclusion and more importantly, it gives us a way to see
through when the plasmon resonance happens if the inclusions (NPs), such as gold (Au) NPs and Silver (Ag) NPs, have negative real part in the parameters.
We show that when the electric permittivity of the NP is such that a contrast parameter becomes the eigenvalue of an
integral operator, plasmon resonance happen. We then analyze the heat generation and transfer process by
strictly deduction. Finally, a typical nanostructure of concentric nanoshell is considered and we show exactly
how is the plasmon resonance excited by incident light.

The organization of this paper is as follows. In section 2, we introduce some preliminary works concerning the
electromagnetic field model. We present the far field expansion of the electromagnetic field perturbed by the
nanoparticles. In the case that the nanoparticles are sphere shaped, we show exactly how the plasmon resonance
happen when Drude model is applied. In section 3, the heat generation and transferring model is considered. The
heat transferring is explicit shown for sphere nanoparticle. Finally, in section 4, we consider a special
concentric nanoshell structure. The strict mathematical proof is given concerning the plasmon resonance. To
our knowledge, it is the first time to show exactly the plasmon resonance in this concentric nanoshell structure.

\section{Plasmon resonance model}
Let $D$ be the nanoparticle situated in $\RR^3$ with
$\mathcal{C}^{1, \eta}$ boundary for some $\eta>0$, and let
$(\epsilon_0,\mu_0)$ be the pair
 of electromagnetic parameters (permittivity and permeability) of the matrix ($\RR^3 \setminus \overline{D}$)
 and $(\epsilon_1, \mu_1)$ be that of the nanoparticle. Then the permittivity and permeability distributions are given by
\begin{equation*}
\epsilon=\epsilon_0 \chi ( \RR^3 \setminus \overline{D}) +
\epsilon_1 \chi(D) \quad \mbox{and} \quad \mu =\mu_0 \chi ( \RR^3
\setminus \overline{D}) + \mu_1 \chi(D),
\end{equation*}
 where $\chi$ denotes the
characteristic function. In the sequel, we set $k_1=\omega
\sqrt{\epsilon_1 \mu_1}$ and $k_0=\omega \sqrt{\epsilon_0 \mu_0}$.
Suppose the nanoparticle $D$ is illuminated by a given incident plane wave with the electric and magnetic fields $(\bE^{in}, \bH^{in})$,
 which is the solution to the Maxwell equations
\begin{equation*}
 \ds\ \left \{
 \begin{array}{ll}
\ds\nabla\times{\bE^{in}} = i\omega\mu_0{\bH^{in}} \quad &\mbox{in } \mathbb{R}^3,\\
\ds\nabla\times{\bH^{in}} = -i\omega\epsilon_0 {\bE^{in}}\quad &\mbox{in
} \mathbb{R}^3
 \end{array}
 \right .
 \end{equation*}
 where $i=\sqrt{-1}$.
We can treat the nanoparticle $D$ as a scatterer in front of the incident wave. The total fields
denoted by $(\bE,\bH)$, is the solution to the following Maxwell
equations: \beq\label{eq:maxwell} \left\{
\begin{array}{ll}
 \nabla \times \bE = i \omega \mu \bH &  \mbox{in} \quad \RR^3\setminus \p D, \\
 \nabla  \times \bH= - i \omega \epsilon \bE & \mbox{in} \quad \RR^3\setminus \p D, \\
 {[}{\bm \nu} \times \bE]= [{\bm \nu} \times \bH] = 0 & \mbox{on} \quad \p D,
\end{array}
\right. \eeq subject to the Silver-M\"{u}ller radiation condition:
 $$\lim_{|\Bx|\rightarrow\infty} |\Bx| (\sqrt{\mu} (\bH- \bH^{in}) \times\hat{\Bx}-\sqrt{\epsilon} (\bE-\bE^{in}))=0,$$
 where $\hat{\bx} = \bx/|\bx|$. Here, $[{\bm \nu} \times{\bE}]$ and $[{\bm \nu} \times{\bH}]$ denote the jump of ${\bm \nu} \times{\bE}$
and ${\bm \nu} \times{\bH}$ along $\p D$, namely,
 $$[{\bm \nu} \times{\bE}]=({\bm \nu} \times \bE)\bigr|^+_{\p D}-({\bm \nu} \times \bE)\bigr|^-_{\p D},
 \quad [{\bm \nu} \times{\bH}]=({\bm \nu} \times \bH)\bigr|^+_{\p D}-({\bm \nu} \times \bH)\bigr|^-_{\p D}$$
where $\cdot|^+_{\p D}$ means the limit to the outside of $\p D$ and $\cdot|^-_{\p D}$ means the limit to the inside of $\p D$.

In what follows, we recall the analytic solution to \eqnref{eq:maxwell}.
Firstly, for  $k>0$, the fundamental solution $\Gamma_k$ to the Helmholtz
operator $(\Delta+k^2)$ in $\RR^3$ is \beq\label{Gk} \ds\Gamma_k
(\bx) = -\frac{e^{ik|\bx|}}{4 \pi |\bx|}. \eeq
To proceed, we introduce some vector functional space and some important integrals
. Let  $\nabla_{\p
D}\cdot$ denote the surface divergence. Denote by $L_T^2(\p D):=\{\bm{\varphi}\in {L^2(\p D)}^3, \bm\nu\cdot \bm{\varphi}=0\}$ the tangential vector space. We introduce the function
spaces
\begin{align*}
TH({\rm div}, \p D):&=\Bigr\{ \bm{\varphi} \in L_T^2(\partial D):
\nabla_{\partial D}\cdot \bm{\varphi} \in L^2(\partial D) \Bigr\},\\
TH({\rm curl}, \p D):&=\Bigr\{ \bm{\varphi} \in L_T^2(\partial D):
\nabla_{\partial D}\cdot (\bm{\varphi}\times \bm{\nu}) \in L^2(\partial D) \Bigr\},
\end{align*}
equipped with the norms
\begin{align*}
&\|\bm{\varphi}\|_{TH({\rm div}, \p D)}=\|\bm{\varphi}\|_{L^2(\p D)}+\|\nabla_{\p D}\cdot\bm{\varphi}\|_{L^2(\p D)}, \\
&\|\bm{\varphi}\|_{TH({\rm curl}, \p D)}=\|\bm{\varphi}\|_{L^2(\p D)}+\|\nabla_{\p D}\cdot(\bm{\varphi}\times\bm\nu)\|_{L^2(\p D)}.
\end{align*}
For a density $\bm{\phi} \in TH({\rm div}, \p D)$, we define the
single layer potential associated with the fundamental solutions
$\Gamma_k$ given in \eqnref{Gk} by
$$
\ds\Scal_D^{k}[\bm{\phi}](\Bx) := \int_{\p D} \Gamma_k(\Bx-\By)
\bm{\phi}(\By) d s(\By), \quad \Bx \in \mathbb{R}^3.
$$
For a scalar density contained in $L^2(\p D)$, the single layer
potential is defined by the same way. We also define boundary
integral operators
\begin{align*}
\ds\mathcal{L}_D^k[\bm{\phi}] (\Bx)& := \bigr(\bm{\nu} \times
\bigr(k^2 \Scal_D^k[
\bm{\phi}] + \nabla \Scal_D^k[\nabla_{\p D}\cdot \bm{\phi}]\bigr) \Bigr)(\Bx),\\
\ds\mathcal{M}_D^k[\bm{\phi}](\Bx) &:= \mbox{p.v.} \int_{\p D}
\bm{\nu}(\Bx) \times \Bigr( \nabla_{\Bx} \times \bigr(
\Gamma_k(\Bx-\By) \bm{\phi} (\By)\bigr)  \Bigr) d s(\By),
\quad \Bx \in \p D.
\end{align*}
There admits the following jump formula on the boundary of $D$
\beq\label{eq:jmp}
\bm{\nu}\times \nabla \times \Scal_D^k[\bm{\varphi}]\Big|_{\pm} = (\mp \frac{I}{2} + \Mcal_D^k)[\bm{\varphi}].
\eeq
 Then the solution to
\eqnref{eq:maxwell} can be represented as the following
 \beq
\label{represent} \ds\bE (\Bx)= \left \{
 \begin{array}{ll}
\ds \bE^i(\Bx) + \mu_0 \nabla \times \Scal_D^{k_0}
[\bm{\phi}](\Bx) +
\nabla\times\nabla\times\Scal_D^{k_0} [\bm{\psi}](\Bx) ,\quad &\Bx \in \mathbb{R}^3 \setminus \overline{D},\\
\nm \ds\mu_1 \nabla \times \Scal_D^{k_1} [\bm{\phi}](\Bx) +
\nabla\times\nabla\times\Scal_D^{k_1} [\bm{\psi}](\Bx) ,\quad &\Bx
\in D,
 \end{array}
 \right .
 \eeq
and
$$ \bH(\Bx) = -\frac{i}{\omega \mu}\bigr(\nabla \times \bE\bigr)(\Bx),\quad \Bx \in \RR^3\setminus \p D,
$$
where the pair $(\bm{\phi}, \bm{\psi}) \in TH({\rm div}, \p D)
\times TH({\rm div}, \p D)$ is the unique solution to \beq
\label{phi_psi}
\begin{bmatrix}
\Mcal_1 &
\ds\mathcal{L}_D^{k_1} - \mathcal{L}_D^{k_0} \\
\ds\mathcal{L}_D^{k_1} - \mathcal{L}_D^{k_0} & \Mcal_2
\end{bmatrix}
\begin{bmatrix}
\bm{\phi} \\ \bm{\psi}
\end{bmatrix}
= \left.\begin{bmatrix}
 \bm{\nu}\times \bE^{in}\\
i \omega \bm{\nu} \times \bH^{in}
\end{bmatrix}\right|_{\p D}.
\eeq
where
$$\Mcal_1:=\ds\frac{\mu_1+\mu_0}{2}I + \mu_1 \mathcal{M}_D^{k_1} -\mu_0
\mathcal{M}_D^{k_0},\quad \Mcal_2:=\ds \left(
\frac{k_1^2}{2 \mu_1} + \frac{k_0^2}{2 \mu_0}\right)I +
\frac{k_1^2}{\mu_1}\mathcal{M}_D^{k_1} -
\frac{k_0^2}{\mu_0}\mathcal{M}_D^{k_0}.$$
 Denote by $\Gl_{\epsilon}$ and $\Gl_{\mu}$ the
electric permittivity and magnetic permeability contrasts: \beq
\Gl_{\epsilon}=\frac{\epsilon_1+\epsilon_0}{2(\epsilon_1-\epsilon_0)}
\quad \mbox{and} \quad
\Gl_{\mu}=\frac{\mu_1+\mu_0}{2(\mu_1-\mu_0)}. \eeq
We point out that if $\Gl_{\epsilon}$ and $\Gl_{\mu}$ are greater than $1/2$,
or the permittivity and permeability are all positive numbers, then the invertibility of the system of equations (\ref{phi_psi})
on $TH({\rm div}, \p D) \times TH({\rm div}, \p D)$ was proved in
\cite{T}.

\subsection{Drude model}
In what follows, we consider the physical process on the plasmon resonance when the nanoparticle $D$ is
illuminated by the incident wave. According to the Drude model the electric permittivity of the nanoparticle
$D$ behaviors as a function of the angular frequency $\omega$, or exactly
(see e.g. \cite{SC10})
\beq\label{eq:Drude}
\epsilon_1=\epsilon(\omega)=\epsilon_0(1-\frac{\omega_p^2}{\omega(\omega+i\tau)})
\eeq
where $\omega_p$ is the plasma frequency of the bulk material and $\tau$ is the width of the resonance.
Let $D=\delta B+ \bz$, where $B$ is a $\mathcal{C}^{1,\eta}$ domain containing the origin. For a scalar density $\phi \in
L^2(\partial B)$, define the
well-known Neumann-Poincar\'e operator by
\begin{equation} \label{defk}
\mathcal{K}_B^*[\phi](\bx) := \mbox{p.v.}\quad \int_{\partial B} \frac {\partial
\Gamma}{\partial {\nu(\bx)}}(\bx-\by) \phi (\by)\,ds(\By),
 \end{equation}
where $\partial/\partial \nu$ denotes the normal derivative and p.v. denotes the Cauchy principle value.
Denote by $\bm{G}(\bx,\bz)$ the matrix valued function (Dyadic Green function)
$$
\bm{G}(\bx,\bz)=\epsilon_0(\Gamma_{k_0}(\bx-\bz)\bm{I} +  \frac{1}{k_0^2}D_{\bx}^2\Gamma_{k_0}(\bx-\bz))
$$
then there holds
the following far field expansion for the electric field
\begin{thm}[Theorem 3.8 in \cite{ADKLMZ}]\label{th:far_expan}
Define the polarization tensors
\beq
\bm{M}^e:=\int_{\p B} \tilde{\by} (\Gl_{\epsilon}I-\Kcal_B^*)^{-1}[\bm\nu] ds(\tilde{\by}) \quad \mbox{and}\quad
\bm{M}^h:=\int_{\p B} \tilde{\by} (\Gl_{\mu}I-\Kcal_B^*)^{-1}[\bm\nu] ds(\tilde{\by})
\eeq
then there holds the following far field expansion
\beq
\bE(\bx)-\bE^{in}(\bx) =-\delta^3\omega^2\mu_0\bm{G}(\bx,\bz)\bm{M}^e\bE^{in}(\bz)-\delta^3\frac{i\omega\mu_0}{\epsilon_0}\nabla\times\bm{G}(\bx,\bz)\bm{M}^h\bH^{in}(\bz) + O(\delta^4).
\eeq
\end{thm}
Since the eigenvalue of Neumann-Poincar\'e operator $\Kcal_B^*$ lies in the span $(-1/2,1/2]$ (cf. \cite{book2}), polarization tensors
$\bm{M}^e$ and $\bm{M}^h$ are well-defined if the material parameters are all positive. However, if the parameters
of the NP are not positive numbers then polarization tensors may not be well-defined since $\Gl_{\epsilon}I-\Kcal_B^*$
and $\Gl_{\mu}I-\Kcal_B^*$ may not be invertible in this case. It is also shown in \cite{ADKLMZ} that when the electric
and magnetic properties of the NPs meet that their real parts make $\Gl_{\epsilon}$ or $\Gl_{\mu}$ be the eigenvalue of
$\Kcal_B^*$, the plasmon resonance happen. We shall see when the parameters obey the Drude and the angular frequency of
the incident wave is well chosen, the plasmon resonance is excited.

\subsection{Sphere nanoparticles}
We shall consider that the nanoparticle $D$ is a sphere shaped inclusion. Suppose the radius of
NP is $r_{NP}$. The electric permittivity of $D$ obeys
the Drude model \eqnref{eq:Drude}. We mention that the wavelength, denoted by $\Gl$, of the incident wave is much larger than the
radius of the NP, i.e., $\Gl>>r_{NP}$. Suppose the incident electric field is uniformly distributed, $\bE^{in}=\bE_0$,
where $\bE_0$ is a constant vector representing the amplitude of the incoming light. Due to the symmetric property
of $D$, we suppose the electric potentials, denoted by $u$, inside and outside $D$ have the form
$$
u=\left\{
\begin{array}{ll}
\bE_0\cdot \bx + \frac{\bE_1\cdot \bx}{|\bx|^3} & \bx \in \RR^3\setminus D,\\
\bE_2\cdot \bx & \bx\in D.
\end{array}
\right.
$$
The scattering part $\bE_1\cdot \bx/|\bx|^3$ decays fast as the light travels far away.
By introducing the spherical harmonic functions $Y_1^m(\theta,\varphi)$, $m=-1, 0, 1$, we
actually have
$$
\bE_0\cdot \bx= r\sum_{m=-1}^1 a_{0m} Y_1^m, \quad \frac{\bE_1\cdot \bx}{|\bx|^3}=r^{-2} \sum_{m=-1}^1 a_{1m} Y_2^m,
\quad \bE_2\cdot \bx= r\sum_{m=-1}^1 a_{2m} Y_1^m
$$
where the coefficients $a_{jm}, j=0, 1, 2$ can be determined by transmission conditions, namely
\begin{align*}
r_{NP}\sum_{m=-1}^1 a_{0m} Y_1^m+r_{NP}^{-2} \sum_{m=-1}^1 a_{1m} Y_2^m=r_{NP}\sum_{m=-1}^1 a_{2m} Y_1^m,\\
\epsilon_0(\sum_{m=-1}^1 a_{0m} Y_1^m-2r_{NP}^{-3} \sum_{m=-1}^1 a_{1m} Y_2^m)=\epsilon_1\sum_{m=-1}^1 a_{2m} Y_1^m.
\end{align*}
By solving the above equations we have
$$
a_{2m}=\frac{3\epsilon_0}{2\epsilon_0+\epsilon_1}a_{0m}, \quad \mbox{and} \quad
a_{1m}=\frac{\epsilon_0-\epsilon_1}{2\epsilon_0+\epsilon_1}r_{NP}^3a_{0m}.
$$
Thus there holds the following relation
$$
\bE_2=\frac{3\epsilon_0}{2\epsilon_0+\epsilon_1}\bE_0, \quad \mbox{and}
\quad \bE_1=\frac{\epsilon_0-\epsilon_1}{2\epsilon_0+\epsilon_1}r_{NP}^3\bE_0.
$$
The energy of the NP, denoted by $ \mathcal{E}$ after the incident light ($\bE_0$) illumination reads
\beq\label{eq:engy}
\mathcal{E}=\int_{D} |\bE_2|^2d\bx=\Big|\frac{3\epsilon_0}{2\epsilon_0+\epsilon_1}\Big|^2\int_D |\bE_0|^2d\bx.
\eeq
In what follows, we show how does the plasmon resonance happen when the electric permittivity of NP obeys the Drude model.
Mathematically, we give the definition on when does the plasmon resonance happen.
\begin{defn}\label{def:1}
Let $\epsilon_1$ satisfy \eqnref{eq:Drude}. We call the plasmon resonance happen if the
energy of the NP induced by the illumination of the incident light blows up in the following way
$$
\lim_{\tau\rightarrow 0} \tau \mathcal{E}=\infty.
$$
\end{defn}
Based on the Definition \ref{def:1} and \eqnref{eq:engy} we can easily get the following result
\begin{thm}
Let $\Gl$ and $v$ be the wave length and speed of the incident light, respectively. If $\Gl=2\pi v\frac{\sqrt{3}}{\omega_p}$,
then the plasmon resonance happen.
\end{thm}
\pf The angular frequency $\omega$ is given by
$$\omega=\frac{2\pi v}{\Gl}=\frac{\sqrt{3}}{3}\omega_p.$$
A direct calculation by using \eqnref{eq:Drude} gives
$$
\epsilon_1=\Big(1-\frac{\omega_p^2}{\omega^2+\tau^2}\Big)\epsilon_0+i\frac{\omega_p^2\tau}{\omega(\omega^2+\tau^2)}\epsilon_0
=-2\epsilon_0+ i\frac{3\sqrt{3}\tau}{\omega_p}\epsilon_0+O(\tau^2).
$$
Then by \eqnref{eq:engy} the energy can be estimated by
$$
\mathcal{E}=\Big|\frac{3\epsilon_0}{i\frac{3\sqrt{3}}{\omega_p}\tau\epsilon_0+O(\tau^2)}\Big|^2\int_D |\bE_0|^2d\bx
=\frac{\sqrt{3}}{3}\omega_p^2\tau^{-2} \int_D |\bE_0|^2d\bx+O(\tau^{-1}).
$$
We thus have
$$
\lim_{\tau\rightarrow 0} \tau \mathcal{E}=\infty
$$
and plasmon resonance happens.
\qed

We mention that in the lossless Drude mode case ($\tau=0$), the permittivity $\epsilon_1$ turns to such that $\Gl_{\epsilon}$
is the eigenvalue of the Neumman-Poincar\'e operator $\Kcal_D^*$. To explain this, we first present a lemma which
gives the eigenvalue of $\Kcal_D^*$ with respect to eigenfunctions $Y_1^m$ when $D$ is a sphere.
\begin{lem}(cf. \cite{ACKLM2})\label{le:le1}
Let $B_0$ be a ball with radius $r_0$ then we have
\begin{align*}
&\Kcal_{B_0}^*[Y_1^m](\bx)=\frac{1}{6}Y_1^m, \quad |\bx|=r=r_0, \quad m=-1,0,1.
\end{align*}
\end{lem}
We thus come to the conclusion since in the lossless Drude mode $\epsilon_1=-2\epsilon_0$, and
$$
\Gl_{\epsilon}=\frac{\epsilon_1+\epsilon_0}{2(\epsilon_1-\epsilon_0)}=\frac{1}{6}
$$
which is exactly the eigenvalue of $\Kcal_D^*$ when the incident light is uniformly distributed.

\section{Heat generation and transfer}
In this section, we consider the heat generated by NPs under plasmon resonance. Heat transfer in
a system with NPs is described by the usual heat transfer equation
\beq\label{eq:heatmodel}
\rho(\bx)c(\bx)\frac{\p T(\bx,t)}{\p t}=\nabla (\Gs(\bx)\nabla T(\bx,t))+ Q(\bx,t)
\eeq
where $T(\bx,t)$ is the temperature, $\rho(\bx)$, $c(\bx)$ and $\Gs(\bx)$ are the mass density,
specific heat, and thermal conductivity, respectively. $Q(\bx,t)$ is the heat intensity which
represents an energy source coming from light dissipation in NPs. It is shown in \cite{GZSRLK2006} that
$Q(\bx,t)$ relates to the electric field by
$$
Q(\bx,t)=\frac{1}{8\pi}\omega \Im m(\epsilon_1)|\bE_2|^2
$$
where $\Im m(\epsilon_1)$ means the imaginary part of the electric permittivity $\epsilon_1$,
given by Drude model and $\bE_2$ has been calculated in the last section. Thus we get
$$
Q(\bx,t)=\frac{1}{8\pi}\omega \Im m(\epsilon_1)\Big|\frac{3\epsilon_0}{2\epsilon_0+\epsilon_1}\Big|^2|\bE_0|^2.
$$
Heat transfers through the NP requires quite few time. We thus consider the steady state of this process.
Denote by $\Gs_0$ and $\Gs_{NP}$ the thermal conductivities outside and insider the NP, respectively.
Then in the steady state, we have
$$
\left\{
\begin{array}{ll}
\Gs_{NP}\Delta T+ Q= 0 & \bx\in D ,\\
\Gs_0\Delta T=0 & \bx\in \RR^3\setminus D.
\end{array}
\right.
$$
Since the heat source is generated from the center of the NP, it spreads
uniformly in every direction due to the uniform thermal diffusion properties of
the NP and the matrix. It also decays to nothing as $|\bx|\rightarrow 0$.
We suppose the temperature inside and outside the NP has the form
$$
T=\left\{
\begin{array}{ll}
A-\frac{Q|\bx|^2}{6\Gs_{NP}} & \bx\in D,\\
B/|\bx| & \bx\in \RR^3\setminus D.
\end{array}
\right.
$$
By using the transmission conditions we thus have
$$
A=\frac{2\Gs_0+\Gs_{NP}}{6\Gs_0\Gs_{NP}}Qr_{NP}^2, \quad B=\frac{r_{NP}^3Q}{3\Gs_0}.
$$
Denote by $V_{NP}$ the volume of the NP then temperature distribution outside the NP is given by
$$
T=\frac{r_{NP}^3Q}{3\Gs_0r}=\frac{V_{NP}Q}{4\pi\Gs_0}\frac{1}{r}, \quad r=|\bx|>r_{NP}
$$
which is in accordance with \cite{GZSRLK2006}. It is seen that the heating effect
becomes especially strong under the plasmon resonance conditions since $Q$ is greatly
increased under plasmon resonance.

\section{Concentric nanoshell}
\begin{figure}[!h]
   \begin{center}
{\includegraphics[width=2.8in]{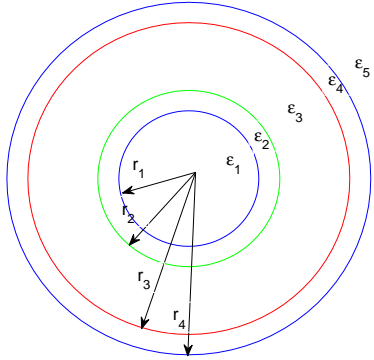}}
   \end{center}
    \caption{Schematic of a concentric nanoshells. $\epsilon_1$ and $\epsilon_3$ are composed of
    dielectric core. $\epsilon_2$ and $\epsilon_4$ are composed of metallic shell. $\epsilon_5$
    is the embedding medium.
       }\label{fig:1}
\end{figure}

In this section, we investigate the plasmonic properties of a concentric nanoshell.
The plasmon hybridization model has been used to explain
the properties of nanoshell, a tunable plasmonic nanoparticle
consisting of a dielectric (silica) core and a metallic (Au or
Ag) shell. Fig. \ref{fig:1} gives an example of concentric nanoshell with two layers of metal
NPs.
This  concentric  nanoshell  consists  of
alternating  layers  of  dielectric  and  metal,  essentially  a
nanoshell enclosed within another nanoshell, inspiring its
alternative  name  of  ¡°nanomatryushka¡± (cf. \cite{RH2004}). We shall show
how does the plasmon resonance happen in this kind of structure. We can also
see how sensitive does the inner and outer radius of the shell layer influence the
plasmon resonance.

From the analysis in the last sections, we should only consider the system in the influence
of the electric field. Suppose once again that the incident wave is uniformly distributed
with amplitude $\bE_0$. Define
\beq
\label{eq:aj}
A_1:=\{r\leq r_1\}, \quad A_j:=\{r_{j}<r\leq r_{j+1}\}, \quad  j=1,2,3, \quad A_5:=\{r>r_4\}.
\eeq
By the symmetric properties of the concentric nanoshell, we suppose
the total electric potential $u$ has the form
$$
u=\left\{
\begin{array}{ll}
\bE_1\cdot \bx= r\sum_{m=-1}^{1}a_{1m}Y_1^m, & \bx\in A_1,\\
\bE_{j1}\cdot \bx+ \bE_{j2}\cdot\bx/|\bx|^3= r\sum_{m=-1}^{1}a_{jm}Y_1^m + r^{-2}\sum_{m=-1}^{1}b_{jm}Y_1^m, & \bx\in A_j,\quad j=2,3,4 ,\\
\bE_0\cdot \bx + \bE_5\cdot \bx/|\bx|^3= r\sum_{m=-1}^{1}a_{0m}Y_1^m+ r^{-2}\sum_{m=-1}^{1}b_{5m}Y_1^m, & \bx\in A_5
\end{array}
\right.
$$
The transmission conditions on the interface $\{ r=r_j\}$, $j=1,2,3,4$ are given by
 \begin{align*}
 &a_{j+1 m} r_j + b_{j+1 m}r_j^{-2} = a_{jm} r_{j} + b_{jm}r_{j}^{-2} , \\
 &\epsilon_{j+1} \left( a_{j+1 m}  - 2b_{j+1 m}r_j^{-3}
\right) = \epsilon_j \left( a_{jm}  - 2b_jr_{j}^{-3}\right)
 \end{align*}
 where we set $b_{1m}=0$ and $a_{5m}=a_{0m}$, $m=-1,0,1$. By setting
  \beq
 \Gl_{j}=\frac{2\epsilon_{j+1}+\epsilon_{j}}{\epsilon_{j+1}-\epsilon_{j}}, \quad j=1,2,3,4
 \eeq
 and some basic arrangements to the equations given by transmission conditions, we get
\begin{align*}
&\Gl_1(a_{2m}-a_{1m})-\sum_{j=2}^4 (a_{j+1 m}-a_{jm})=-a_{0m} \\
&-2(a_{2m}-a_{1m})\Big(\frac{r_1}{r_2}\Big)^3-\Gl_2(a_{3m}-a_{2m})+\sum_{j=3}^4 (a_{j+1 m}-a_{jm})=a_{0m} \\
&2\sum_{j=1}^2(a_{j+1 m}-a_{jm})\Big(\frac{r_j}{r_3}\Big)^3+\Gl_3(a_{4m}-a_{3m})-(a_{5m}-a_{4m})=-a_{0m}\\
&-2\sum_{j=1}^3(a_{j+1 m}-a_{jm})\Big(\frac{r_j}{r_4}\Big)^3-\Gl_4(a_{5m}-a_{4m})=a_{0m}
\end{align*}
Define the matrix $P$ and $\Upsilon_{n}$ by
 \beq
 P:= \begin{bmatrix}
 \Gl_1 & -1 & -1 & -1 \\
-2(r_{1}/r_2)^3 & -\Gl_{2} & 1 & 1 \\
2(r_{1}/r_3)^3 & 2(r_{2}/r_3)^3 & \Gl_3 & -1 \\
-2(r_1/r_4)^3 & -2(r_2/r_4)^3 & -2(r_{3}/r_4)^3 & -\Gl_4
 \end{bmatrix}
 \quad
 \Upsilon:= \begin{bmatrix}
 r_1^3 & 0 & 0 & 0 \\
0 & r_{2}^3 & 0 & 0 \\
0 & 0 & r_3^3 & 0 \\
0 & 0 & 0 & r_4^3
 \end{bmatrix} .
\eeq
Then there holds
\beq\label{eq:b}
\Bb_m = a_{0m}\Xi \Upsilon P^{-1} \itbf{e}
\eeq
where $\Bb_m:=(b_{2m},b_{3m},b_{4m},b_{5m})^T$, $\itbf{e}:=(1,-1,1,-1)^T$ and
$$
\Xi=
\begin{bmatrix}
1 & 0 & 0 & 0 \\
1 & 1 & 0 & 0 \\
1 & 1  & 1 & 0 \\
1 & 1  & 1 & 1
\end{bmatrix}.
$$
We also have
\beq\label{eq:a}
\itbf{a}_m = a_{0m}(-\Xi^TP^{-1}\itbf{e} +\itbf{e}_4)
\eeq
where $\itbf{a}_m:=(a_{1m},a_{2m},a_{3m},a_{4m})^T$ and $\itbf{e}_4:=(1,1,1,1)^T$.
What we concern is the energy generated in the metal shells. By the analysis before,
we have in mind that in plasmon resonance mode, the energy blows up if the imaginary part
of material property of the metal particles goes to zero.
To simplify the analysis, in what follows we suppose that the material properties of the
dielectric core are fixed and $\epsilon_1=\epsilon_3=\epsilon_5$
Let $\epsilon_{s}=\epsilon_2=\epsilon_4$. Then we have
$$\Gl_1=\Gl_3=1-\Gl_2=1-\Gl_4=\frac{2\epsilon_s+\epsilon_1}{\epsilon_s-\epsilon_1}.$$
Under these assumptions, we have $P=\Gl_1 I-K$, where $I$ is the identity matrix and $K$ has the form
 \beq
 K= \begin{bmatrix}
0 & 1 & 1 & 1 \\
2(r_{1}/r_2)^3 & 1 & -1 & -1 \\
-2(r_{1}/r_3)^3 & -2(r_{2}/r_3)^3 & 0 & 1 \\
2(r_1/r_4)^3 & 2(r_2/r_4)^3 & 2(r_{3}/r_4)^3 & 1
 \end{bmatrix}
\eeq
Judging from \eqnref{eq:b} and \eqnref{eq:a}, we set the determinant of $P$, or $\Gl_1I-K$, to be zero and
we have the equation
\begin{align}
0=&\quad\quad\Gl_1^4 - 2\Gl_1^3+\left(-2r_1^3r_2^{-3}+ 2r_1^3r_3^{-3}- 2r_1^3r_4^{-3}+ 2r_2^3r_4^{-3} - 2r_3^3r_4^{-3}+ 1\right)\Gl_1^2 +\nonumber\\
& \left(2r_1^3r_2^{-3}- 2r_1^3r_3^{-3}+2r_1^3r_4^{-3} +2 r_2^3r_3^{-3}- 2r_2^3r_4^{-3}   + 2r_3^3r_4^{-3}\right)\Gl_1 + 4r_1^3r_2^{-3} r_3^3r_4^{-3}\label{eq:4th}
\end{align}
which is a forth order equation with respect to $\Gl_1$. It is easy to see that the solution are
exactly the eigenvalues of the matrix $K$. We can actually solve the forth order equation explicitly by using
some tools like matlab, etc. It can be seen that all the four roots of the equation are real. However, we shall not
discuss the reason here why all solutions should be real but the reader who is interested may consult any book
concerning the theory of algebra equations (see, e.g., \cite{Ti01}).
We also see in this equation that the solution only depends on the ratios of the radii.
Table \ref{tab:1} lists the solutions to four different kinds of ratios.
\begin{table}[h]
\caption{Solutions to the fourth order equation \eqnref{eq:4th}.}
\begin{center}
\begin{tabular}{ccc|ccc}
\hline\noalign{\smallskip}
 $r_1:r_2:r_3:r_4$ & $\Gl_1$ & $\epsilon_s/\epsilon_1$ & $r_1:r_2:r_3:r_4$ & $\Gl_1$ & $\epsilon_s/\epsilon_1$ \\
 \noalign{\smallskip} \hline \noalign{\smallskip}
 & -0.5915 &  -0.1576                   &            & -0.5013 &  -0.1994\\
$4:5:9:10$ &    1.5915 & -6.3439 & $3:4:7:8$  & 1.5013 &  -5.0151\\
 & -0.8550 &   -0.0508                &            & -0.8237&   -0.0624\\
  &  1.8550 &    -19.6878              &            & 1.8237&  -16.0206\\
\noalign{\smallskip} \hline \noalign{\smallskip}
 & -0.6546 &  -0.1301                   &            & -0.3787 &  -0.2612  \\
 $5:6:11:12$ &   1.6546 & -7.6849  & $3:4:6:8$  & 1.3787 & -3.8288\\
  & -0.8771 &   -0.0427               &            & -0.7702 &  -0.0830\\
   & 1.8771 &  -23.4024                &            & 1.7702 &   -12.0547\\
\noalign{\smallskip} \hline \noalign{\smallskip}
\end{tabular}
\end{center}
\label{tab:1}
\end{table}

The following theorem gives exactly how the plasmon resonance happen in this concentric nanoshell structure.
\begin{thm}
Suppose $\epsilon_1=\epsilon_3=\epsilon_5$
Let $\epsilon_{s}=\epsilon_2=\epsilon_4$. If the angular frequency $\omega$ of the incident light is chosen
such that
$$\lim_{\tau\rightarrow 0}\epsilon_{s}(\omega)=\epsilon^*(\omega)$$ and there exists a vector $\itbf{v}_1^*\in \RR^4$
such that
$$P(\epsilon^*)\itbf{v}_1^*=0, \quad |\itbf{v}_1^*|=1\quad \mbox{and} \quad \itbf{v}_1^*\cdot \itbf{e}\neq 0$$
then the plasmon resonance happen in concentric nanoshell.
\end{thm}
\pf By the assumptions and Drude model we have
$$\epsilon_{s}(\omega)=\epsilon^*+i\tau\frac{\omega_p}{\omega^3}+O(\tau^2).$$
We then have
$$
\Gl_1(\epsilon_s)=\Gl_1(\epsilon^*)+i\tau \frac{\omega_p}{\omega^3}\Big(-(2\epsilon^*+\epsilon_1)+\frac{1}{\epsilon^*-\epsilon_1}\Big)+O(\tau^2).
$$
We also have that $\Gl_1(\epsilon^*)$ is one eigenvalue of $K$, with corresponding eigenvector $\itbf{v}_1^*$. Furthermore,
$$
P(\epsilon_s)=P(\epsilon^*)+i\tau \frac{\omega_p}{\omega^3}\Big(-(2\epsilon^*+\epsilon_1)+\frac{1}{\epsilon^*-\epsilon_1}\Big) I + O(\tau^2)I
$$
 In the following we let $b:=\frac{\omega_p}{\omega^3}\Big(-(2\epsilon^*+\epsilon_1)+\frac{1}{\epsilon^*-\epsilon_1}\Big)$. By the assumptions we have
$$
P(\epsilon_s)\itbf{v}_1^*=(i\tau b+O(\tau^2))\itbf{v}_1^*
$$
and thus
$$
P(\epsilon_s)^{-1}\itbf{v}_1^*=(-i\tau^{-1}b^{-1}+O(1))\itbf{v}_1^*.
$$
Next, we calculate the energy generated in the metal shell. Denote by $\mathcal{E}$ the energy
generated in the metal shell then we have
$$
\mathcal{E}=\int_{A_2}\left|\sum_{m=-1}^1 (ra_{2m}Y_1^m+r^{-2}b_{2m}Y_1^m)\right|^2d\bx
+\int_{A_4}\left|\sum_{m=-1}^1 (ra_{4m}Y_1^m+r^{-2}b_{4m}Y_1^m)\right|^2d\bx
$$
By using the orthogonality of the spherical harmonic functions $Y_1^m$, \eqnref{eq:b} and \eqnref{eq:a} we have
\begin{align*}
\mathcal{E}&=\frac{r_2^3-r_1^3}{3}\sum_{m=-1}^1 |a_{2m}|^2+\frac{r_4^3-r_3^3}{3}\sum_{m=-1}^1 |a_{4m}|^2-\frac{r_2^{-5}-r_1^{-5}}{5}\sum_{m=-1}^1 |b_{2m}|^2-\frac{r_4^{-5}-r_3^{-5}}{5}\sum_{m=-1}^1 |b_{4m}|^2 \\
&=\Big(\frac{r_2^3-r_1^3}{3}|f_1(\epsilon_s)|^2+\frac{r_4^3-r_3^3}{3}|f_2(\epsilon_s)|^2+\frac{r_1^{-5}-r_2^{-5}}{5}|f_3(\epsilon_s)|^2
+\frac{r_3^{-5}-r_4^{-5}}{5}|f_4(\epsilon_s)|^2\Big)\sum_{m=-1}^1 |a_{0m}|^2
\end{align*}
where the four functions which depend on $\epsilon_s$, i.e., $f_1(\epsilon_s)$, $f_2(\epsilon_s)$, $f_3(\epsilon_s)$, $f_3(\epsilon_s)$ have
the form
\begin{align*}
&f_1(\epsilon_s):=-(0,1,1,1)P(\epsilon_s)^{-1}\itbf{e}+1,  \quad f_2(\epsilon_s):=-(0,0,0,1)P(\epsilon_s)^{-1}\itbf{e}+1 & \\
&f_3(\epsilon_s):=(1,0,0,0)\Upsilon P(\epsilon_s)^{-1}\itbf{e}, \quad \quad f_4(\epsilon_s):=(1,1,1,0)\Upsilon P(\epsilon_s)^{-1}\itbf{e}.&
\end{align*}
Define
$$
F=
\begin{bmatrix}
0 & -1  & -1 & -1 \\
0 & 0  & 0 & -1 \\
r_1^3 & 0 & 0 & 0 \\
r_1^3 & r_2^3 & r_3^3 & 0
\end{bmatrix}
$$
then $(f_1(\epsilon_s), f_2(\epsilon_s), f_3(\epsilon_s), f_4(\epsilon_s))^T= FP(\epsilon_s)^{-1}\itbf{e}+(1,1,0,0)^T$.
Since we have $\det(F)=r_1^3(r_2^3-r_3^3)\neq 0$, with out loss of generality we suppose $(1,0,0,0)\Upsilon \itbf{v}_1^*\neq 0$.
Since all the four terms in the expression of $\mathcal{E}$ are positive terms,
we only analyze the term which contains $f_3(\epsilon_s)$ (otherwise, if $(0,1,1,1)\itbf{v}_1^*\neq 0$ then we consider $f_1(\epsilon_s)$ and
so on).
Let $\itbf{v}_1^*$ ,$\itbf{v}_2^*$, $\itbf{v}_3^*$ and $\itbf{v}_4^*$ be the eigenvectors of $K$  and $\Gl_1(\epsilon^*)$, $\Gl_2^*$
, $\Gl_3^*$, $\Gl_4^*$ be the corresponding eigenvalues. Then
$\itbf{e}= \itbf{v}_1^*\itbf{e}\cdot \itbf{v}_1^* + \itbf{v}_2^*\itbf{e}\cdot \itbf{v}_2^* +\itbf{v}_3^*\itbf{e}\cdot \itbf{v}_3^*  +\itbf{v}_4^*\itbf{e}\cdot \itbf{v}_4^* $.
Direct calculations give
\begin{align*}
P(\epsilon_s)^{-1}\itbf{e}&=(\Gl_1(\epsilon_s)-\Gl_1(\epsilon^*))^{-1}\itbf{v}_1^*\itbf{e}\cdot \itbf{v}_1^*+\sum_{j=2}^4 (\Gl_1(\epsilon_s)-\Gl_j^*)^{-1}\itbf{v}_j^*\itbf{e}\cdot \itbf{v}_j^* \\
&= -i\tau^{-1} b^{-1}\itbf{v}_1^*\itbf{e}\cdot \itbf{v}_1^* +O(1).
\end{align*}
Thus we have
$$
\tau |f_3(\epsilon_s)|^2=\tau |-i\tau^{-1} b^{-1}(1,0,0,0)\Upsilon\itbf{v}_1^*\itbf{e}\cdot \itbf{v}_1^* +O(1)|^2
=\tau^{-1}b^{-2}|(1,0,0,0)\Upsilon\itbf{v}_1^*|^2|\itbf{e}\cdot \itbf{v}_1^*|^2+O(1).
$$
By the assumptions $(1,0,0,0)\Upsilon\itbf{v}_1^*\neq 0$ and $\itbf{e}\cdot \itbf{v}_1^*\neq 0$ we finally have
$$
\lim_{\tau \rightarrow 0} \tau \mathcal{E} =\infty
$$
which completes the proof by using Definition \ref{def:1}.
\qed

We make a short remark here. We have seen in Table \ref{tab:1} that there are four
possible frequencies, according to Drude model, of incident waves which can be used
to excite the plasmon resonance. The concentric nanoshell structure, in the mean time,
only contains two metal shells. Thus if the number of metal shells and dielectric cores
is increasing, more different frequencies of incident lights can be absorbed and turned
to a large amount of heat. The photothermal effects can be used for a variety of imaging
and therapeutic applications.


\section{Conclusions}
We studied the heat generation and transferring model in the presence of plasmon resonance
when the NPs are illuminated by incident waves. For a electromagnetic wave illumination, the
first order expansion in terms of polarization tensors was presented, and the key point for
inducing plasmon resonance emerged. We showed strictly how the plasmon resonance happen when
the nanoparticle is sphere shaped and obeys the Drude model. We investigated the heat
generation and transferring for spherical nanoparticle. The photothermal effect is greatly
enhanced under plasmon resonance. For a concentric nanoshell structure, we proved that the
plasmon resonance happen when the frequency of the incident waves is well chosen. Future works
will be focused on the interaction of the nanoparticles under plasmon resonance which is not
only a very important physical problem but also a great mathematical problem.


\begin{thebibliography}{1}
\bibitem{ACKLM2} H. Ammari, G. Ciraolo, H. Kang, H. Lee, G.W. Milton, Spectral theory
of a Neumann-Poincar\'e-type operator and analysis of anomalous localized resonance II, arXiv:1212.5066.

\bibitem{ADKLMZ}
H. Ammari, Y. Deng, H. Kang, H. Lee, P. Millien, and J. Zou, Analysis of plasmon resonances in nano-particles, preprint.

\bibitem{book2} {H. Ammari and H. Kang}, \textsl{Polarization and Moment
Tensors with Applications to Inverse Problems and Effective Medium
Theory}, Applied Mathematical Sciences, Vol. 162, Springer-Verlag,
New York, 2007.

\bibitem{AVV2001}
H. Ammari,  M. Vogelius, and D. Volkov, Asymptotic formulas for
perturbations in the electromagnetic fields due to the presence of
inhomogeneities of small diameter II.
 The full Maxwell equations,
 J. Math. Pures Appl., {80} (2001), 769--814.

\bibitem{BL2002}
M. A. Barral and A. M. Llois, Photothermal Imaging of
Nanometer-Sized Metal Particles among Scatterers,
Science,  297 (2002), 1160--1163.

\bibitem{BLBCL2006}
S. Berciaud, D. Lasne, G. A. Blab, L. Cognet and B. Lounis,
Photothermal Heterodyne Imaging of Individual Metallic
Nanoparticles: Theory versus Experiment, Phys. Rev. B
73 (2006), 045424.


\bibitem{ESMLM1997}
R. Elghanian, J.J. Storhoff, R.C. Mucic, R.L. Letsinger and C.A. Mirkin,
Selective Colorimetric Detection of Polynucleotides Based on the Distance-Dependent Optical Properties of Gold Nanoparticles,
Science, 22 (1997), 1078--1081.

\bibitem{GR07}
A. O. Govorov and H. Richardson, Generating heat with metal nanoparticles, NanoToday, (1) 2 (2007), 20-39.

\bibitem{GZSRLK2006}
A. O. Govorov, W. Zhang, T. Skeini, H. Richardon, J. Lee and N. A. Kotov,
Gold nanoparticle ensembles as heaters and acuators: melting and collective plasmon resonances,
Nanoscale Res. Lett., 1 (2006), 84-90.

\bibitem{GLHJDW2007}
A. M. Gobin, M. H. Lee, N. J. Halas, W. D. James, R. A. Drezek and
West, J. L., Near-Infrared Resonant Nanoshells for
Combined Optical Imaging and Photothermal Cancer
Therapy, Nano Lett., 7 (2007), 1929--1934.



\bibitem{HEQE2006}
X. Huang , I.H. El-Sayed , W. Qian and M.A. El-Sayed,
Cancer Cell Imaging and Photothermal Therapy in the Near-Infrared Region by Using Gold Nanorods,
J. Am. Chem. Soc., (6) 128 (2006), 2115--2120.

\bibitem{HF04}
E. Hutter and J.H. Fendler, Exploitation of localized surface plasmon resonance,
Adv. Mater., {16} (2004), 1685-1706.

\bibitem{JS07}
P.K. Jain, I.H. El-Sayed, and M.A. El-Sayed, Au nanoparticles target cancer,
Nanotoday, {2} (2007), 18-29.

\bibitem{Ka13}
M. Kalteh, Investigating the effect of various nanoparticle and base liquid types on the nanofluids heat and fluid flow in a microchannel,
Applied Mathematical Modelling,  (18-19) 37 (2013), 8600--8609.

\bibitem{MNB1993}
C. B. Murray, D. J. Norris and M. G. Bawendi,
Synthesis and characterization of nearly monodisperse CdE (E = sulfur, selenium, tellurium) semiconductor nanocrystallites,
J. Am. Chem. Soc., (1993) 115, 8706--8715.


\bibitem{Pe99}
J. B. Pendry, Radiative Exchange of Heat Between Nanostructures, J. Phys.: Condens. Matter {11} (1999), 6621-6633.

\bibitem{RH2004}
C. Radloff and N. J. Halas, Plasmonic Properties of Concentric Nanoshells, Nano
letters, (7) 4 (2004), 1323--1327.

\bibitem{SC10}
D. Sarid and W. A. Challener, \textsl{Modern Introduction to Surface Plasmons: Theory, Mathematica Modeling, and Applications},
Cambridge University Press, New York, 2010.

\bibitem{Ti01}
J. P. Tignol, \textsl{Galois' Theory of Algebraic Equations}, World Scientific Publishing, Singapore, 2001.

\bibitem{T} R. H. Torres, {Maxwell's equations and dielectric obstacles
with Lipschitz boundaries}, J. London Math. Soc., (2) 57 (1998),
157-169.

\bibitem{WMTH2009}
B. Wu, S. W. McCue, P. Tillman, J. M. Hill, Single phase limit for melting nanoparticles, 
Applied Mathematical Modelling, (5) 33 (2009), 2349--2367.

\bibitem{ZK2009}
W. Zhao and J. M. Karp,  Tumour Targeting: Nanoantennas
Heat Up, Nat. Mater., 8(2009), 453--454.


\end{thebibliography}
\end{document}